\documentclass[11pt]{article}
\usepackage[latin9]{inputenc}
\usepackage{geometry}
\usepackage{amsmath}
\usepackage{amssymb}
\usepackage{esint}
\usepackage[authoryear]{natbib}

\makeatletter
\usepackage{amsfonts}
\usepackage{graphicx}
\usepackage{caption}
\usepackage{subcaption}
\usepackage{amsfonts}
\usepackage{amsthm}

\theoremstyle{definition}
\newtheorem{theorem}{Theorem}\newtheorem{corollary}{Corollary}\newtheorem{lemma}{Lemma}\theoremstyle{remark}
\newtheorem{example}{Example}\newtheorem{definition}{Definition}\newtheorem{remark}{Remark}

\makeatother

\begin{document}

\title{Invertibility of infinitely divisible continuous-time moving average
processes}

\author{Orimar Sauri\thanks{Financial support from the Center for Research in the Econometric
Analysis of Time Series (grant DNRF78) funded by the Danish National
Research Foundation is gratefully acknowledged. This study was also
partially funded by the Villum Fonden as part of the project number
11745 titled \textit{''Ambit Fields: Probabilistic Properties and
Statistical Inference''}.} \\
 %EndAName
Department of Mathematics and CREATES\\
 Aarhus University\\
 osauri@math.au.dk }

\date{\today }
\maketitle
\begin{abstract}
This paper studies the invertibility property of continuous time moving
average processes driven by a Lévy process. We provide of sufficient
conditions for the recovery of the driving noise. Our assumptions
are specified via the kernel involved and the characteristic triplet
of the background driving Lévy process.
\end{abstract}

\section{Introduction}

In the context of time series, the concept of \textit{invertibility}
of stochastic processes refers to the task of recovering the driving
noise by the observed series. Such a property plays an important role
for the characterization of the notion of \textit{causality}, which
is the principle in where the current state of a given system is not
influenced by its future states. Invertibility and causality are well
understood in the discrete-time framework, in particular, for moving
average processes, necessary and sufficient conditions for invertibility
and causality have been established in terms of its moving average
coefficients. See for instance \citet{Brockwell:1986:TST:17326}.
Motivated by this framework, the main goal of the present paper is
to study the invertibility property of the class of continuous-time\textit{
moving average processes driven by a Lévy process}, that is, the observed
process $\left(X_{t}\right)_{t\in\mathbb{R}}$ admits the spectral
representation
\begin{equation}
X_{t}:=\int_{\mathbb{R}}f\left(t-s\right)dL_{s},\text{ \ \ }t\in\mathbb{R}\text{,}\label{eqn0.1}
\end{equation}
where $f$ is a measurable function, often called \textit{kernel},
and $L$ is a Lévy process. Our main result states that the process
$X$ is invertible, for a certain class of Lévy processes, whenever
the Fourier transform of $f$ does not vanish, which is in essence
the analogous condition to the discrete-time setting. We would like
to emphasized that the class of Lévy processes we consider in our
results does not need to be square integrable. See Section \ref{sec:Invertibility-of-IDCMA}
for more details.

Observe that the process $X$ is infinitely divisible in the sense
of \citet{idprocessesmaejimasatoole} and \citet{BNSauSzo15}. Thus,
in statistical terms, the kernel $f$ models the autocorrelation structure
of $X$ while $L$ describe its distributional properties. Furthermore, $X$
can be used as a flexible model that is able to reproduce many of
the stylized properties found in empirical data such as fat tails
and local Gaussianity (mixed Gaussian distributions). Hence, from
the modeling perspective, invertibility provides a simple way to identify
(in a one-to-one relation) and estimate the law of $X$ by $L$, and
vice versa.

Several authors have investigated the invertibility problem for continuous-time
processes. For instance, \citet{invnoncausacomte} studied the invertibility
and causality of Gaussian Volterra processes, which are those processes
that can be written as in (\ref{eqn0.1}) but we replace $f(t-s)$
by $f(t,s)$ and $L$ by a Brownian motion. Under smoothness assumptions
on the kernel, the authors provided necessary and sufficient conditions
for the invertibility and causality of these type of processes. In
the non-Gaussian case, \citet{CohMaej2011} established the invertibility
property for the the family of fractional Lévy processes in the case
when $L$ is centered and has finite second moment. 

In the stationary framework, \citet{Brockwell20092660} considered
the continuous-time version of the classical ARMA processes. In their
set up, the authors gave necessary and sufficient conditions (which
turned out to be the analogous of those for the classical ARMA) for
the causality and invertibility of this family. Recently, \cite{BasseNielPedRoh17}
studied the solutions of ARMA type stochastic differential equations.
The authors showed that when the solution exists, it can be written
as in (\ref{eqn0.1}) and, under extra regularity conditions, such
a solution is invertible and causal. The previous situations are contained
in our framework.

The present paper is organized as follows. Section 2 introduces the
notation and some background on infinite divisibility, stochastic
integration with respect to Lévy processes, and Orlicz spaces. In
Section 3, we present our main result and we discuss several important
examples. Section 4 concludes.

\section{Preliminaries and basic results}

Throughout this paper $\left(\Omega,\mathcal{F},\left(\mathcal{F}_{t}\right)_{t\in\mathbb{R}},\mathbb{P}\right)$
denotes a filtered probability space satisfying the usual conditions
of right-continuity and completeness. For $p\geq0$, we denote by
$\mathcal{L}^{p}\left(\Omega,\mathcal{F},\mathbb{P}\right)$ the space
of $p$-integrable random variables endowed with the convergence in
$p$-mean for $p>0$ and convergence in probability for the case when
$p=0$.

A two-sided $\mathbb{R}^{d}$-valued Lévy process $\left(L_{t}\right)_{t\in\mathbb{R}}$
on $\left(\Omega,\mathcal{F},\mathbb{P}\right)$ is a stochastic process
taking values in $\mathbb{R}^{d}$ with independent and stationary
increments whose sample paths are almost surely càdlàg. We say that
$\left(L_{t}\right)_{t\in\mathbb{R}}$ is an $\left(\mathcal{F}_{t}\right)$-Lévy
process if for all $t>s,$ $L_{t}-L_{s}$ is $\mathcal{F}_{t}$-measurable
and independent of $\mathcal{F}_{s}$.

By $ID\left(\mathbb{R}^{d}\right)$ we mean the space of infinitely
divisible distributions on $\mathbb{R}^{d}$. Any Lévy process is
infinitely divisible and $L_{1}$ has a Lévy-Khintchine representation, relative to a truncation function $\tau$,
given by 
\[
\log\widehat{\mu}\left(z\right)=i\left\langle z,\gamma_{\tau}\right\rangle -\frac{1}{2}\left\langle z,Bz\right\rangle +\int_{\mathbb{R}^{d}}\left[e^{i\left\langle z,x\right\rangle }-1-i\left\langle \tau\left(x\right),z\right\rangle \right]\nu\left(\mathrm{d}x\right),\text{ \ \ }z\in\mathbb{R}^{n}\text{,}
\]
where $\widehat{\mu}$ is the characteristic function of the law of
$L_{1}$, $\gamma_{\tau}\in\mathbb{R}^{d}$, $B$ is a symmetric nonnegative
definite matrix on $\mathbb{R}^{d\times d}$, and $\nu$ is a Lévy
measure, i.e. $\nu\left(\left\{ 0^{d}\right\} \right)=0$, with $0^{d}$
denoting the origin in $\mathbb{R}^{d},$ and $\int_{\mathbb{R}^{d}}(1\wedge\left\vert x\right\vert ^{2})\nu\left(\mathrm{d}x\right)<\infty.$
Here, we assume that the truncation function $\tau$ is given by $\tau\left(x_{1},\ldots,x_{n}\right)=\left(\frac{x_{i}}{1\vee\left\vert x_{i}\right\vert }\right)_{i=1}^{n},\ \ \left(x_{1},\ldots,x_{n}\right)\in\mathbb{R}^{n}$.

An \emph{infinitely divisible continuous-time moving average} (IDCMA)
process is a stochastic process $\left(X_{t}\right)_{t\in\mathbb{R}}$
on $\left(\Omega,\mathcal{F},\left(\mathcal{F}_{t}\right)_{t\in\mathbb{R}},\mathbb{P}\right)$
given by the following formula
\begin{equation}
X_{t}:=\int_{\mathbb{R}}f\left(t-s\right)dL_{s},\text{ \ \ }t\in\mathbb{R},\label{eqn1.1}
\end{equation}
where $f$ is a deterministic function and $L$ is a Lévy process
with triplet $\left(\gamma_{\tau},B,\nu\right)$. IDCMA process belongs
to the class of \textit{Lévy semistationary process} ($\mathcal{LSS}$)
which are those processes $\left(Y_{t}\right)_{t\in\mathbb{R}}$ which
are described by the following dynamics
\begin{equation}
Y_{t}=\theta+\int_{-\infty}^{t}g\left(t-s\right)\sigma_{s}dL_{s}+\int_{-\infty}^{t}q\left(t-s\right)a_{s}ds,\text{ \ \ }t\in\mathbb{R}\text{,}\label{eqn1.2}
\end{equation}
where $\theta\in\mathbb{R}^{d}$, $L$ is a Lévy process, $g$ and
$q$ are deterministic functions such that $g\left(x\right)=q\left(x\right)=0$
for $x\leq0$, and $\sigma$ and $a$ are adapted càdlàg processes.
For further references to theory and applications of Lévy semistationary
processes, see \citet{RePEc:aah:create:2010-18} and references therein.

\subsection{Stochastic integrals and Orlicz spaces}

In the following, we present a short review of \citet{rajputrosinski}
and \citet{satostochinadditive} concerning the existence of stochastic
integrals of the form $\int_{\mathbb{R}}f(s)dL_{s}$, where $f:\mathbb{R}\rightarrow\mathbb{R}$
is a measurable function and $L$ a Lévy process as well as the connection
of such integral with the so-called Orlicz spaces.

Let $L$ be an $\mathbb{R}^{d}$-valued Lévy process with characteristic
triplet $(\gamma_{\tau},B,\nu)$. The space of simple functions on
$\mathbb{R}$ will be denoted by $\vartheta$. Thus, $f\in\vartheta$
if and only if $f$ can be written as 
\[
f=\sum\limits _{i=1}^{k}a_{i}\mathbf{1}_{(s_{i},t_{i}]}\text{,}
\]
where $s_{i}\leq t_{i}$ and $a_{i}\in\mathbb{R}$ for $i=1,\ldots,k$.
For any $f\in\vartheta$, the integral of $f$ with respect to (w.r.t.
for short) is defined as
\[
\int_{\mathbb{R}}f(s)dL_{s}:=\sum\limits _{i=1}^{k}a_{i}(L_{t_{i}}-L_{s_{i}})\text{.}
\]
We will say that $f$ is $L$-integrable if there exists a sequence
$(f_{n})_{n\geq1}\subseteq\vartheta$, such that $f_{n}\rightarrow f$
almost everywhere and that the sequence $\int_{\mathbb{R}}f_{n}(s)dL_{s}$
has a limit in probability and we write
\[
\int_{\mathbb{R}}f(s)dL_{s}:=\mathbb{P}\text{-}\lim_{n\rightarrow\infty}\int_{\mathbb{R}}f_{n}(s)dL_{s},
\]
In \citet{rajputrosinski}, c.f. \citet{satostochinadditive}, it
has been shown that $f$ is $L$-integrable and $\int_{\mathbb{R}}f(s)dL_{s}\in\mathcal{L}^{p}\left(\Omega,\mathcal{F},\mathbb{P}\right)$
if and only if $\int_{\mathbb{R}}\Phi_{p}^{(\gamma_{\tau},B,\nu)}(f(s))ds<\infty,$
where 
\begin{equation}
\Phi_{p}^{(\gamma_{\tau},B,\nu)}(u):=V(u)+tr(B)u^{2}+\int_{\mathbb{R}^{d}}[\left\Vert ux\right\Vert ^{2}\mathbf{1}_{\left\Vert ux\right\Vert \leq1}+\left\Vert ux\right\Vert ^{p}\mathbf{1}_{\left\Vert ux\right\Vert >1}]\nu(dx),\,\,\,u\in\mathbb{R},\label{eqn1.4}
\end{equation}
with
\[
V(u):=\left\vert \gamma_{\tau}u+\int_{\mathbb{R}^{d}}\left[\tau\left(ux\right)-u\tau\left(x\right)\right]\nu(dx)\right\vert ,\,\,\,u\in\mathbb{R}\text{.}
\]
Observe that for $p>0$, $\Phi_{p}^{(\gamma_{\tau},B,\nu)}$ is well
defined if and only if $\int_{\left\Vert x\right\Vert >1}\left\Vert x\right\Vert ^{p}\nu(dx)<\infty$.
For the rest of this paper the space of $L$-integrable functions
will be denoted by
\[
\mathbb{L}_{\Phi_{p}^{(\gamma_{\tau},B,\nu)}}:=\{f:(\mathbb{R},\mathcal{B}(\mathbb{R}))\rightarrow(\mathbb{R},\mathcal{B}(\mathbb{R})):\int_{\mathbb{R}}\Phi_{p}^{(\gamma_{\tau},B,\nu)}(\left|f(s)\right|)ds<\infty\}.
\]
In general, $\mathbb{L}_{\Phi_{p}^{(\gamma_{\tau},B,\nu)}}$ is a
complete linear metric space in which $\vartheta$ is dense, but it
is not necessarily a Banach space. However, under certain conditions
on $\Phi_{p}^{(\gamma_{\tau},B,\nu)}$, $\mathbb{L}_{\Phi_{p}^{(\gamma_{\tau},B,\nu)}}$
becomes equivalent to an \emph{Orlicz Space}, which is a certain type
of Banach space. Hence, we now present some properties of such spaces.
We refer the reader to \citet{Orliczrao} for more details.

A mapping $\Psi:\mathbb{R}\rightarrow\left[0,\infty\right]$ is said
to be a \emph{Young function} if it is even, convex with $\Psi(s)=0$
if and only if $s=0$, and such that $\lim_{s\rightarrow\infty}\Psi(s)=+\infty$.
Given a Young function $\Psi$, the mapping 
\begin{equation}
\overline{\Psi}\left(x\right):=\sup_{y\geq0}\left\{ \left\vert x\right\vert y-\Psi\left(y\right)\right\} \text{, \ \ }x\in\mathbb{R}\text{.}\label{eqn1.6-1}
\end{equation}
define a new Young function and it is termed as its complementary
function. We say that a function $\Psi$ fulfills the $\Delta_{2}$-condition
if $\Psi\left(2x\right)\leq K\Psi\left(x\right)$ for some $K>0$.
For a given Young function satisfying the $\Delta_{2}$-condition
let
\[
\mathcal{L}_{\Psi}:=\left\{ f:(\mathbb{R},\mathcal{B}(\mathbb{R}))\rightarrow(\mathbb{R},\mathcal{B}(\mathbb{R})):\int_{\mathbb{R}}\Psi\left(\left\vert f\left(s\right)\right\vert \right)ds<\infty\right\} .
\]
We have that in this framework, $\mathcal{L}_{\Psi}$ is a separable
Banach space equipped with Luxemburg norm
\begin{equation}
\left\Vert f\right\Vert _{\Psi}:=\inf\left\{ a>0:\int_{\mathbb{R}}\Psi\left(a^{-1}\left\vert f\left(s\right)\right\vert \right)\mathrm{d}s\leq1\right\} ,\label{eqn1.6.0-1}
\end{equation}
when equivalent functions are identified almost everywhere. $\mathcal{L}_{\Psi}$
is known as the Orlicz space associated to $\Psi$. 
By $\mathcal{S}(\mathbb{R})$
we mean the space of \emph{test functions of rapidly decaying}, i.e.
$\phi\in\mathcal{S}(\mathbb{R})$ if it is infinitely continuously
differentiable and for any $n\geq1$ and $m\geq0$, the mapping $x\mapsto\phi^{(m)}(x)x^{n}$
is bounded on $\mathbb{R}$, where $\phi^{(m)}$ denotes the derivative
of order $m$ of $\phi$. The space of \emph{tempered distributions},
which we denote by $\mathcal{S}^{\prime}(\mathbb{R})$, is the topological
dual of $\mathcal{S}(\mathbb{R})$. For more details on the theory
of tempered distributions we refer to \cite{DuisKolk10}. Fix $\Psi$
a non-trivial Young function, i.e. $\Psi(x)\neq+\infty$, $x>0$,  satisfying
the $\Delta_{2}$-condition. We have the following connections between
Orlicz spaces and the the space of tempered distributions:
\begin{enumerate}
\item Let $f\in\mathcal{L}_{\Psi}$, then $f$ is locally integrable and
by Jensen's inequality, for any $n\geq1$
\[
\Psi(\int_{\mathbb{R}}\left|\frac{f(s)}{(1+\left|s\right|)^{n}}\right|ds)\leq c_{n}\int_{\mathbb{R}}\Psi(\left|f(s)\right|)ds<\infty.
\]
The latter, according to \cite{DuisKolk10}, p. 189, gives us that
$\mathcal{L}_{\Psi}\subseteq\mathcal{S}^{\prime}(\mathbb{R})$.
\item If $f\in\mathcal{L}_{\Psi},g\in\mathcal{L}_{\overline{\Psi}}$. Then
for any $t\in\mathbb{R}$
\[
\int_{\mathbb{R}}\left|f(t-s)g(s)\right|ds\leq2\left\Vert f\right\Vert _{\Psi}\left\Vert g\right\Vert _{\overline{\Psi}}.
\]
For a proof see \cite{Orliczrao}, p. 58.
\item By the previous point, if $f\in\mathcal{L}_{\Psi},g\in\mathcal{L}_{\overline{\Psi}}$,
we get that for any $n\geq1$
\[
\int_{\mathbb{R}}\left|\frac{f\ast g(s)}{(1+\left|s\right|)^{n}}\right|ds<\infty,
\]
which means that the induced distribution by $f\ast g$ belongs to
$\mathcal{S}^{\prime}(\mathbb{R})$.
\end{enumerate}
The next result identify $\mathcal{L}_{\Psi}^{\prime}$, the dual
of $\mathcal{L}_{\Psi}$, 

\begin{theorem}[\cite{Orliczrao}, p. 105.]The dual of $\mathcal{L}_{\Psi}$
is isometrically isomorphic to $\mathcal{L}_{\overline{\Psi}}$, where
$\overline{\Psi}$ is as in (\ref{eqn1.6-1}). More precisely, for
any $T\in\mathcal{L}_{\Psi}^{\prime}$ there exists a unique $g\in\mathcal{L}_{\overline{\Psi}}$,
such that 
\[
T(f)=\int_{\mathbb{R}}f(s)g(s)ds\text{, \ \ }f\in\mathcal{L}_{\Psi}.
\]
\end{theorem}

Recall that in a Banach space $\left(\mathcal{X},\left\Vert \cdot\right\Vert _{\mathcal{X}}\right),$
a collection $F=\left(f_{\alpha}\right)_{\alpha\in\Lambda}$ is said
to be dense if $\overline{F}=\mathcal{X}$ under the norm $\left\Vert \cdot\right\Vert _{\mathcal{X}}$.
From the previous theorem and the Hahn-Banach Theorem we get:

\begin{corollary} \label{totalityorlicz}A collection $F=\left(f_{\alpha}\right)_{\alpha\in\Lambda}\subset\mathcal{L}_{\Psi}$
is dense in $\mathcal{L}_{\Psi}$ if and only if 
\[
\int_{\mathbb{R}}f_{\alpha}(s)g(s)ds=0,\text{\ \ \ }\forall\text{ }\alpha\in\Lambda\text{,}
\]
with $g\in\mathcal{L}_{\overline{\Psi}}$, implies that $g\equiv0$,
almost everywhere. \end{corollary}

Turning back to the stochastic integral, fix $p\geq0$ and suppose
that $\Phi_{p}^{(\gamma_{\tau},B,\nu)}$ is comparable to a Young
function, that is, there are $c,C>0$ and a Young function $\Psi$,
such that 
\begin{equation}
c\Psi(x)\leq\Phi_{p}^{(\gamma_{\tau},B,\nu)}(x)\leq C\Psi(x),\,\,\,x\geq0.\label{youngcondition}
\end{equation}
Since $\Phi_{p}^{(\gamma_{\tau},B,\nu)}$ satisfies the $\Delta_{2}$-condition
(\citet{rajputrosinski}), we conclude that in this case $\mathcal{L}_{\Psi}$
is a Banach space equivalent to $\mathbb{L}_{\Phi_{p}^{(\gamma_{\tau},B,\nu)}}$.

\begin{remark} We observe the following:
\begin{enumerate}
\item Although the Lévy processes under consideration are $\mathbb{R}^{d}$-valued,
the space $\left(\mathcal{L}_{\Psi},\left\Vert \cdot\right\Vert _{\Psi}\right)$
contains only real-valued functions.
\item From \citet{Kaminska1997}, an Orlicz space $\left(\mathcal{L}_{\Psi},\left\Vert \cdot\right\Vert _{\Psi}\right)$
is isometric to some Hilbert space if and only if $\Psi(x)=kx^{2}$
for come $k>0$. Therefore, $\mathbb{L}_{\Phi_{p}^{(\gamma_{\tau},B,\nu)}}$
is comparable to a Hilbert space if and only if $L$ is centered and
square integrable.
\end{enumerate}
\end{remark} 

The following properties of the stochastic integral defined above
will be useful for the rest of the paper, see \citet{rajputrosinski}
for a proof:

\begin{theorem} \label{continuitystochint}Let $\left(L_{t}\right)_{t\in\mathbb{R}}$
be a Lévy process with triplet $\left(\gamma_{\tau},B,\nu\right)$
and suppose that (\ref{youngcondition}) holds for some $p\geq0$.
Then 
\begin{enumerate}
\item The mapping $\left(f\in\mathcal{L}_{\Psi}\right)\mapsto\left(\int_{\mathbb{R}}f(s)dL_{s}\in\mathcal{L}^{p}\left(\Omega,\mathcal{F},\mathbb{P}\right)\right)$
is continuous, i.e. if $\left\Vert f_{n}-f\right\Vert _{\Psi}\rightarrow0$
, then $\int_{\mathbb{R}}f_{n}(s)dL_{s}\rightarrow\int_{\mathbb{R}}f(s)dL_{s}$
in $\mathcal{L}^{p}\left(\Omega,\mathcal{F},\mathbb{P}\right)$;
\item If $L$ is symmetric, then $\left(f\in\mathcal{L}_{\Psi}\right)\mapsto\left(\int_{\mathbb{R}}f(s)dL_{s}\in\mathcal{L}^{p}\left(\Omega,\mathcal{F},\mathbb{P}\right)\right)$
is an isomorphism between $\mathcal{L}_{\Psi}$ and $\mathcal{L}^{p}\left(\Omega,\mathcal{F},\mathbb{P}\right)$,
that is, if $\int_{\mathbb{R}}f_{n}(s)dL_{s}\rightarrow\int_{\mathbb{R}}f(s)dL_{s}$
in probability, then $\left\Vert f_{n}-f\right\Vert _{\Psi}\rightarrow0$.
Moreover
\[
\overline{\mathrm{span}}\{L_{t}-L_{s}:s\leq t\}=\{\int_{\mathbb{R}}f(s)dL_{s}:f\in\mathcal{L}_{\Psi}\},
\]
where the closure is taken on $\mathcal{L}^{p}\left(\Omega,\mathcal{F},\mathbb{P}\right)$.
\end{enumerate}
\end{theorem}

\section{Invertibility of IDCMA processes\label{sec:Invertibility-of-IDCMA}}

In this section we present the main result of this paper. Let us start
by recalling the notions of invertibility and causality in the time
series framework. Let $\left(X_{t}\right)_{t\in\mathbb{Z}}$ be a
discrete-time moving average process, i.e. 
\[
X_{t}=\sum\limits _{j\in\mathbb{Z}}\theta_{j}\varepsilon_{t-j}=\Theta\left(B\right)\varepsilon_{t},\text{ \ \ }t\in\mathbb{Z}\text{,}
\]
where the process $\left(\varepsilon_{t}\right)_{t\in\mathbb{Z}}$
is a mean zero weak stationary white noise, $\sum\limits _{j\in\mathbb{Z}}\left\vert \theta_{j}\right\vert <\infty$,
$B$ is the lag operator and 
\[
\Theta\left(z\right)=\sum\limits _{j\in\mathbb{Z}}\psi_{j}z^{j},\text{ \ 	\ }z\in\mathbb{C},\left\vert z\right\vert <1.
\]
Observe that if $\Theta^{-1}$ admits a power series expansion, then
almost surely
\begin{equation}
\varepsilon_{t}=\Theta^{-1}\left(B\right)X_{t}=\sum\limits _{j\in\mathbb{Z}}\pi_{j}X_{t-j},\text{ \ \ }t\in\mathbb{Z}\text{.}\label{eq:invertitimeseries}
\end{equation}
Thus, $\varepsilon_{t}\in\overline{\mathrm{span}}\left\{ X_{s}\right\} _{s\in\mathbb{Z}}$
for any $t\in\mathbb{Z}$, where the closure is taken in $\mathcal{L}^{2}\left(\Omega,\mathcal{F},\mathbb{P}\right)$,
or in other words $X$ is \emph{invertible}. A necessary and sufficient
condition for the power series expansion of $\Theta^{-1}$ is that
$\Theta$ does not vanish in the unitary circle. Thus, (\ref{eq:invertitimeseries})
holds if and only if $\Theta\left(z\right)\neq0$ for all $z\in\mathbb{C},$
with $\left\vert z\right\vert \leq1$. Observe that the latter is
equivalently to 
\[
0\neq\sum\limits _{j\in\mathbb{Z}}\theta_{j}e^{-ij\omega}=\Theta\left(e^{-\omega}\right)=:\widehat{\Theta}\left(\omega\right),\,\,\,\forall\,\left\vert \omega\right\vert \leq\pi.
\]
Note that $\widehat{\Theta}$ is the discrete Fourier transform of
the moving average coefficients $\left(\theta_{j}\right)_{j\in\mathbb{Z}}$.
Hence, the Fourier transform of $\left(\theta_{j}\right)_{j\in\mathbb{Z}}$
does not vanish if and only if (\ref{eq:invertitimeseries}) is satisfied.
These ideas can be extended to characterize the situations in which
$\varepsilon_{t}\in\overline{\mathrm{span}}\left\{ X_{s}\right\} _{s\in\mathbb{Z}}$
for any $t\in\mathbb{Z}$, see \cite{Brockwell:1986:TST:17326} for
more details.

Invertibility itself does not tell us anything about adaptability
of the process. For instance, if $X$ follows an autorregresive dynamics,
i.e. 
\[
X_{t}=\theta X_{t-1}+\varepsilon_{t},\text{ \ \ }t\in\mathbb{Z}\text{,}
\]
then $X$ has stationary a solution if and only if $\theta\neq1$.
In particular, if $\left\vert \theta\right\vert <1$
\begin{equation}
X_{t}=\sum\limits _{j\geq0}\theta^{j}\varepsilon_{t-j},\text{ \ \ }t\in\mathbb{Z}\text{,}\label{eqn2.1}
\end{equation}
and if $\left\vert \theta\right\vert >1$
\begin{equation}
X_{t}=-\sum\limits _{j\geq0}\theta^{-j}\varepsilon_{t+j},\text{ \ \ }t\in\mathbb{Z}\text{.}\label{eqn2.2}
\end{equation}
Note that in (\ref{eqn2.1}), $X$ only depends on the past innovations
of $\varepsilon$ contrary to that in (\ref{eqn2.2}), in which $X$
is expressed in terms of the future innovations of $\varepsilon$.
When $X$ admits a representation as in (\ref{eqn2.1}), it is called
causal and for the case of (\ref{eqn2.2}) it is called non-causal.
However, it is obvious that $\varepsilon$ only depends on the past
innovations of $X$, i.e. $\varepsilon$ admits a causal representation.
This property is usually called\textit{ invertibility in the causal
sense}.

In analogy with the discrete-time framework, we introduce the notion
of invertibility for an IDCMA.

\begin{definition} \label{defininvcaus}Let $X$ be as in (\ref{eqn1.1}).
$X$ is said to be \emph{invertible} on $\mathcal{L}^{p}\left(\Omega,\mathcal{F},\mathbb{P}\right)$
for some $p\geq0$, if $L_{t}-L_{s}\in\overline{\mathrm{span}}\left\{ X_{u}\right\} _{u\in\mathbb{R}}$
for any $t>s$, where the closure is taken in $\mathcal{L}^{p}\left(\Omega,\mathcal{F},\mathbb{P}\right)$.
In the same context, we are going to say that $X$ is\textit{ invertible
in the causal sense} if $L_{t}-L_{s}\in\overline{\mathrm{span}}\left\{ X_{u}\right\} _{u\leq t}$
for any $t>s$.\end{definition}

A natural question appears, as in the discrete-time case, is $\widehat{f}\neq0$
a sufficient (necessary) condition for the invertibility of an IDCMA?
In the case when $\mathbb{L}_{\Phi_{p}^{(\gamma_{\tau},B,\nu)}}$
is equivalent to an Orlicz space, the answer is affirmative as the
following theorem shows.

\begin{theorem} \label{maintheorem}Let $\left(L_{t}\right)_{t\in\mathbb{R}}$
be a Lévy process with characteristic triplet $\left(\gamma,B,\nu\right)$
and suppose that for some $p\geq0$, there is a Young function $\Psi$
satisfying (\ref{youngcondition}). If $f\in\mathcal{L}_{\Psi}\cap\mathcal{L}^{1}\left(dx\right)$
has non-vanishing Fourier transform, then
\begin{equation}
\overline{\mathrm{span}}\left\{ X_{u}\right\} _{u\in\mathbb{R}}=\overline{\mathrm{span}}\left\{ L_{t}-L_{s}:s\leq t\right\} ,\text{ \ \ in }\mathcal{L}^{p}\left(\Omega,\mathcal{F},\mathbb{P}\right)\text{.}\label{eqn1.8}
\end{equation}
\end{theorem}

\begin{remark} In a very informal way, Theorem \ref{maintheorem}
says that for every $t\geq s$ there exists a measurable function
$g_{t,s}$ such that $L_{t}-L_{s}=\int_{\mathbb{R}}g_{t,s}(r)dX_{r}$.
However, since $X$ is not in general a semimartingale, such integral
may not be well defined. \end{remark}

Before presenting the proof of this theorem, we discuss several important
examples.

\begin{example}[Symmetric and integrable L\'evy processes]\label{symmetriccentered}Suppose
that $L$ is a symmetric Lévy process with $\mathbb{E}(\left\Vert L_{1}\right\Vert )<\infty$.
Thus, if $L$ has characteristic triplet $\left(\gamma,B,\nu\right)$,
then 
\[
\Phi_{1}(u):=tr(B)u^{2}+\int_{\mathbb{R}^{d}}(\left\Vert ux\right\Vert ^{2}\wedge\left\Vert ux\right\Vert )\nu(dx),\,\,\,u\in\mathbb{R}.
\]
From the proof of Theorem 3.3 in \cite{BasseRosinski13}, we have
that the mapping 
\[
\Psi(u):=tr(B)u^{2}+\int_{\mathbb{R}^{d}}[\left\Vert ux\right\Vert ^{2}\mathbf{1}_{\left\Vert ux\right\Vert \leq1}+2(\left\Vert ux\right\Vert -1)\mathbf{1}_{\left\Vert ux\right\Vert >1}]\nu(dx),
\]
 is convex and such that 
\[
\Psi(u)/2\leq\Phi_{1}(u)\leq\Psi(u),\,\,\,u\in\mathbb{R}.
\]
Therefore $L$ satisfies the assumptions of Theorem \ref{maintheorem}
(i.e. $\Psi$ is a Young function) if $B\neq0$ or $B=0$ and as $u\rightarrow\infty$
\[
\int_{\mathbb{R}^{d}}(\left\Vert ux\right\Vert ^{2}\wedge\left\Vert ux\right\Vert )\nu(dx)\rightarrow+\infty.
\]

\end{example}

\begin{example}[Ornstein-Uhlenbeck processes]\label{exampleOU}Let
$L$ be a Lévy process with characteristic triplet $\left(\gamma_{\tau},B,\nu\right)$
and put 
\[
f\left(s\right):=e^{-s}\mathbf{1}_{\left\{ s\geq0\right\} }\text{, \ \ }s\in\mathbb{R}\text{.}
\]
Then $X$, the resulting IDCMA process, is the classic OU process
driven by $L$. It is well known that $f\in\mathbb{L}_{\Phi_{0}^{(\gamma_{\tau},B,\nu)}}$
if and only if $\int_{\left\vert x\right\vert >1}\log\left(\left\vert x\right\vert \right)\nu\left(dx\right)<\infty.$
Moreover, since $\widehat{f}$, the Fourier transform of $f$, never
vanishes, we conclude that $f$ satisfies the assumptions of Theorem
\ref{maintheorem}. Furthermore, due to the Langevin equation, it
follows that $X$ is in fact invertible in the causal sense. Now,
if consider instead the process
\[
X_{t}^{\prime}:=\int_{t}^{\infty}e^{-\left(s-t\right)}dL_{s}\text{, \ \ }t\in\mathbb{R},
\]
we get that $X^{\prime}$ is not adapted but well defined provided
that $\int_{\left\vert x\right\vert >1}\log\left(\left\vert x\right\vert \right)\nu\left(dx\right)<\infty.$
Nevertheless, it is easy to check that $X$ fulfills a sort of Langevin
equation, that is, almost surely,
\[
\int_{s}^{t}X_{r}dr=L_{t}-L_{s}+X_{t}-X_{s}\text{, \ \ }t\geq s\text{.}
\]
Hence, we deduce that $X$ is invertible in the causal sense. Observe
that the Langevin equation holds in a pathwise sense, so for the invertibility
of OU-type processes, the condition (\ref{youngcondition}) is superfluous.

\end{example}

\begin{example}[$\mathcal{LSS}$ with a Gamma kernel]\label{examplegamma}Denote
by $L$ a Lévy process with characteristic triplet $\left(\gamma_{\tau},B,\nu\right)$.
Let $\alpha>-1$ and consider
\begin{equation}
f\left(s\right):=e^{-\lambda}s^{\alpha}\mathbf{1}_{\left\{ s>0\right\} }\text{, \ \ }s\in\mathbb{R}\text{.}\label{eqn1.3.15}
\end{equation}
It has been shown in \citet{gammaandreas}, c.f. \citet{PedSau15},
that $f\in\mathbb{L}_{\Phi_{0}^{(\gamma_{\tau},B,\nu)}}$ if and only
if the following two conditions are satisfied:
\begin{enumerate}
\item $\int_{\left\vert x\right\vert >1}\log\left(\left\vert x\right\vert \right)\nu\left(dx\right)<\infty$,
\item One of the following conditions holds:
\begin{enumerate}
\item $\alpha>-1/2;$
\item $\alpha=-1/2$, $B=0$ and $\int_{\left\vert x\right\vert \leq1}\left\vert x\right\vert ^{2}\left\vert \log\left(\left\vert x\right\vert \right)\right\vert \nu\left(dx\right)<\infty;$
\item $\alpha\in\left(-1,-1/2\right)$, $B=0$ and $\int_{\left\vert x\right\vert \leq1}\left\vert x\right\vert ^{-1/\alpha}\nu\left(dx\right)<\infty.$ 
\end{enumerate}
\end{enumerate}
On the other hand, if $p>0$, we claim that $f\in\mathbb{L}_{\Phi_{p}^{(\gamma_{\tau},B,\nu)}}\cap\mathbb{L}_{\Phi_{0}^{(\gamma_{\tau},B,\nu)}}$
if and only if $\alpha p>-1$ and $\int_{\left\vert x\right\vert >1}\left\Vert x\right\Vert ^{p}\nu\left(dx\right)<\infty$.
Indeed, we first observe that there are $c,C>0$ such that 
\[
c\phi_{\alpha,\lambda/2}(s)\leq f(s)\leq C\phi_{\alpha,\lambda}\left(s\right),\,\,\,s>0,
\]
where
\[
\phi_{\alpha,\lambda}\left(s\right):=\left\{ \begin{array}{cc}
s^{\alpha}\mathbf{1}_{\left\{ 0<s\leq1\right\} }+e^{-\lambda s}\mathbf{1}_{\left\{ s>1\right\} } & \text{for }-1/2<\alpha<0;\\
e^{-\lambda s}\mathbf{1}_{\left\{ s\geq0\right\} } & \text{for }\alpha\geq0.
\end{array}\right.
\]
Hence $f\in\mathbb{L}_{\Phi_{p}^{(\gamma_{\tau},B,\nu)}}\cap\mathbb{L}_{\Phi_{0}^{(\gamma_{\tau},B,\nu)}}$
if and only if $\phi_{\alpha,\lambda}\in\mathbb{L}_{\Phi_{p}^{(\gamma_{\tau},B,\nu)}}\cap\mathbb{L}_{\Phi_{0}^{(\gamma_{\tau},B,\nu)}}$.
Our claim then follows by noting that for $\alpha\geq0$ 
\[
\int_{0}^{\infty}\int_{\mathbb{R}^{d}}\left\Vert \phi_{\alpha,\lambda}(s)x\right\Vert ^{p}\mathbf{1}_{\left\Vert \phi_{\alpha,\lambda}(s)x\right\Vert >1}\nu\left(dx\right)ds=\frac{1}{\lambda p}\int_{\left\Vert x\right\Vert >1}\left\Vert x\right\Vert ^{p}(1-\left\Vert x\right\Vert ^{-1})\nu\left(dx\right),
\]
while for $\alpha p>-1$
\begin{align*}
\int_{0}^{\infty}\int_{\mathbb{R}^{d}}\left\Vert \phi_{\alpha,\lambda}(s)x\right\Vert ^{p}\mathbf{1}_{\left\Vert \phi_{\alpha,\lambda}(s)x\right\Vert >1}\nu\left(dx\right)ds & =\frac{1}{p\alpha+1}\int_{\left\vert x\right\vert >1}\left\Vert x\right\Vert ^{p}\nu\left(dx\right)\\
 & +\frac{1}{p\alpha+1}\int_{\left\vert x\right\vert \leq1}\left\Vert x\right\Vert ^{-1/\alpha}\nu\left(dx\right)\\
 & +\frac{1}{\lambda p}\int_{\left\Vert x\right\Vert >e}\left\Vert x\right\Vert ^{p}(\left\Vert x\right\Vert ^{-1}-e)\nu\left(dx\right).
\end{align*}
In this case $X$, the associated IDCMA process, is called \emph{Lévy
semistationary process with a gamma kernel}. See \citet{PedSau15}
for more properties on this process. Note that the Fourier transform
of $f$ is given by 
\[
\hat{f}\left(\xi\right)=\frac{\Gamma\left(\alpha+1\right)}{\sqrt{2\pi}}\frac{1}{\left(\lambda+i\xi\right)^{\alpha+1}}\text{, \ \ }\xi\in\mathbb{R}\text{.}
\]
Hence, under the framework of Theorem \ref{maintheorem}, $X$ is
invertible. Furthermore, it is possible to show that if $\int_{\left\vert x\right\vert >1}\left\Vert x\right\Vert \nu\left(dx\right)<\infty$,
then for any $-1<\alpha<0$, almost surely
\begin{equation}
\int_{0}^{\infty}X_{t-u}\mu\left(\mathrm{d}u\right)=k_{\alpha}\int_{-\infty}^{t}e^{-\lambda\left(t-s\right)}\mathrm{d}L_{s},\text{ \ \ for any }t\in\mathbb{R}\text{,}\label{eqn5.6-1}
\end{equation}
where $\mu\left(du\right):=e^{-\lambda u}u^{-\alpha-1}\left(u\right)\mathbf{1}_{\left\{ u\geq0\right\} }\mathrm{d}u$
and $k_{\alpha}>0$. This relation actually shows that $X$ is invertible
in the causal sense provided that $\int_{\left\vert x\right\vert >1}\left\Vert x\right\Vert \nu\left(dx\right)<\infty$.
As final remark we would like to mention that equation (\ref{eqn5.6-1})
was originally proved in \citet{RePEc:aah:create:2010-18} for the
case when $L$ is a subordinator.\end{example}

\begin{example}[$CARMA(p,q)$]\label{examplecarma}The Lévy driven
$CARMA(p,q)$ (continuous-time auto-regressive moving average process)
with parameters $p>q$, constitutes the generalization of the classical
ARMA models in time series to the continuout-time framework. They
were introduced in \citet{Brockwell20092660} as the stationary process
given by $X_{t}=\mathbf{b}^{\prime}Y_{t}$ where $Y$ follows the
following SDE
\[
dY_{t}=AY_{t}dt+e_{p}dL_{t}\text{,}
\]
where $L$ is a real-valued Lévy process with characteristic triplet
$\left(\gamma,B,\nu\right)$, $\mathbf{b}=\left(b_{0},\ldots,b_{p-1}\right)^{\prime}$,
$\mathbf{e}_{p}=\left(0,0,\cdots,1\right)^{\prime}$ and
\[
A=\begin{bmatrix}0 & 1 & 0 & \cdots & 0\\
0 & 0 & 1 & \cdots & 0\\
\vdots & \vdots & \vdots & \ddots & \vdots\\
0 & 0 & 0 & \cdots & 1\\
-a_{p} & a_{p-1} & a_{p-2} & \cdots & -a_{1}
\end{bmatrix}.
\]
where $a_{1},\ldots,a_{p},b_{0},\ldots,b_{p-1}$ are such that $b_{q}\neq0$
and $b_{j}=0$ for $j>q$. The authors showed that $X$ can be written
as an IDCMA
\[
X_{t}=\int_{\mathbb{R}}g\left(t-s\right)dL_{s}\text{, \ \ }t\in\mathbb{R}\text{,}
\]
with 
\[
g\left(s\right)=b^{\prime}e^{As}\mathbf{e}_{p}\mathbf{1}_{\left\{ s>0\right\} }\text{,}
\]
provided that $\int_{\left\vert x\right\vert >1}\log\left(\left\vert x\right\vert \right)\nu\left(dx\right)<\infty$
and the roots of the polynomial $a\left(\lambda\right)=a_{p}+a_{p-1}\lambda+\cdots+a_{1}\lambda^{p-1}+\lambda^{p},$
$\lambda\in\mathbb{C}$, have strictly negative real part. Since in
this case
\[
\widehat{g}\left(\xi\right)=\frac{b\left(-i\xi\right)}{a\left(-i\xi\right)}\text{, \ \ }\xi\in\mathbb{R}\text{,}
\]
with $b\left(\lambda\right)=b_{0}+b_{1}\lambda+\cdots+b_{p-1}\lambda^{p-1},$
$\lambda\in\mathbb{C}$, we conclude that the kernel of a $CARMA(p,q)$
satisfies the assumptions of Theorem \ref{maintheorem} if the roots
of the polynomial $b$ have non-vanishing real part, i.e. if $b\left(\lambda^{\ast}\right)=0$
then $Re\lambda^{\ast}\neq0$, and $a$ and $b$ have no common roots.
Observe that this condition coincides with the Assumption 1 in \citet{ferrazzano2013}.
For generalizations on the CARMA equation introduced before we refer
to \cite{BasseNielPedRoh17}. \end{example}

The proof of Theorem \ref{maintheorem} in mainly based on the following
lemma.

\begin{lemma} \label{lemma1}Let $\left(L_{t}\right)_{t\in\mathbb{R}}$
be a Lévy process with characteristic triplet $\left(\gamma_{\tau},B,\nu\right)$
and $\Psi$ as in Theorem \ref{maintheorem}. Let $\left(f_{\alpha}\right)_{\alpha\in\Lambda}\subset\mathcal{L}_{\Psi}$.
If $\mathbf{1}_{\left(s,t\right]}\in\mathrm{\overline{span}}\left(f_{\alpha}\right)_{\alpha\in\Lambda}$
under $\left\Vert \cdot\right\Vert _{\Psi}$ for $s\leq t$, then
$L_{t}-L_{s}\in\mathrm{\overline{span}}\left(\int_{\mathbb{R}}f_{\alpha}\left(s\right)dL_{s}\right)_{\alpha\in\Lambda}$
in $\mathcal{L}^{p}\left(\Omega,\mathcal{F},\mathbb{P}\right)$. \end{lemma}

\begin{proof} If $\mathbf{1}_{\left(s,t\right]}\in\mathrm{\overline{span}}\left(f_{\alpha}\right)_{\alpha\in\Lambda}$
under $\left\Vert \cdot\right\Vert _{\Psi}$ for $s\leq t$, then
there exist $\mathbf{\theta}^{n}:=\left(\theta_{i}^{n}\right)_{i=1}^{n}\in\mathbb{R}^{n}$
and $\alpha^{n}:=\left(\alpha_{i}^{n}\right)_{i=1}^{n}\subset\Lambda$
with $n\in\mathbb{N}$, such that $\left\Vert \sum_{i=1}^{n}\theta_{i}^{n}f_{\alpha_{i}^{n}}-\mathbf{1}_{\left(s,t\right]}\right\Vert _{\Psi}\rightarrow0$.
Therefore, from Theorem \ref{continuitystochint}, for some $p\geq0$,
\[
\int_{\mathbb{R}}\sum_{i=1}^{n}\theta_{i}^{n}f_{\alpha_{i}^{n}}(r)dL_{r}\rightarrow L_{t}-L_{s},\,\,\,\,in\,\mathcal{L}^{p}\left(\Omega,\mathcal{F},\mathbb{P}\right),
\]
which is enough. \end{proof}

\begin{proof}[Proof of Theorem \ref{maintheorem}] Obviously $\overline{\mathrm{span}}\left\{ X_{u}\right\} _{u\in\mathbb{R}}\subseteq\overline{\mathrm{span}}\left\{ L_{t}-L_{s}:s\leq t\right\} $
so we only need to show the opposite contention. Recall that under
our assumptions, for some $p\geq0$, $\mathbb{L}_{\Phi_{p}^{(\gamma_{\tau},B,\nu)}}$
is equivalent to the Orlicz space $(\mathcal{L}_{\Psi},\left\Vert \cdot\right\Vert _{\Psi})$.
Thus, from Lemma \ref{lemma1}, we only need to check
that for every $u>s$, $\mathbf{1}_{\left(s,u\right]}\in\mathrm{\overline{span}}\left(f(t-\cdot)\right)_{t\in\mathbb{R}}$
under $\left\Vert \cdot\right\Vert _{\Psi}$. We will prove something
stronger, namely 
\begin{equation}
\mathrm{\overline{span}}\left\{ f\left(t-\cdot\right)\right\} _{t\in\mathbb{R}}=\mathcal{L}_{\Psi}.\label{spanproofmth}
\end{equation}
To do this we will apply Corollary \ref{totalityorlicz}.
Thus, let $g\in\mathcal{L}_{\overline{\Psi}}$ in such a way that
\[
\int_{\mathbb{R}}f\left(t-s\right)g\left(s\right)ds=0\text{, \ for all }t\in\mathbb{R}\text{.}
\]
From Section 2 we know that the functions $f,g$ and $f\ast g$ induce distributions
on $\mathcal{S}^{\prime}(\mathbb{R})$. Thus, their distributional
Fourier transforms are well defined. Denote by $sp(\hat{g})\text{ and }sp(\hat{f})$
the (distributional) support of the Fourier transforms of $g$ and $f$, respectively.
Since $f\in\mathcal{L}_{\Psi}\cap\mathcal{L}^{1}\left(dx\right)$,
we can apply Lemma 5 in \citet{troungvan}, c.f. \cite{MR1466349},
to get that 
\[
sp(\hat{g})\subseteq sp(\hat{f})^{c}=\emptyset,
\]
This implies immediately that $g\equiv0$ almost everywhere, which
according to Corollary \ref{totalityorlicz}, gives (\ref{spanproofmth}).
\end{proof}

\begin{remark} Observe that the reasoning in the previous proof holds
for any integrable function $f:\mathbb{R}^{d}\rightarrow\mathbb{R}$,
having non-vanishing Fourier transform. Therefore, Theorem \ref{maintheorem}
is also applicable for random fields of the form (\ref{eqn1.1}).
\end{remark}

\section{Conclusions}

This paper studied the invertibility of continuous-time moving averages
processes driven by a Lévy processes. We show that driving noise can
be recovered by direct observations of the process. To do this we
assumed that the Fourier transform of the kernel never vanishes and
we imposed a regularity condition on the characteristic triplet of
the background driving Lévy process.

\section*{Acknowledgement}

The author gratefully acknowledges to Ole E. Barndorff-Nielsen and
Benedykt Szozda for helpful comments on a previous version of this
work. 

\bibliographystyle{chicago}
\bibliography{bibIDCMA}

\begin{thebibliography}{}

\bibitem[\protect\citeauthoryear{Bang}{Bang}{1997}]{MR1466349}
Bang, H.~H. (1997).
\newblock Spectrum of functions in {O}rlicz spaces.
\newblock {\em J. Math. Sci. Univ. Tokyo\/}~{\em 4\/}(2), 341--349.

\bibitem[\protect\citeauthoryear{Barndorff-Nielsen, Sauri, and
  Szozda}{Barndorff-Nielsen et~al.}{2015}]{BNSauSzo15}
Barndorff-Nielsen, O., O.~Sauri, and B.~Szozda (2015).
\newblock Selfdecomposable fields.
\newblock {\em To appear in Journal of Theoretical Probability\/}.

\bibitem[\protect\citeauthoryear{Barndorff-Nielsen, Benth, and
  Veraart}{Barndorff-Nielsen et~al.}{2013}]{RePEc:aah:create:2010-18}
Barndorff-Nielsen, O.~E., F.~E. Benth, and A.~Veraart (2013).
\newblock Modelling energy spot prices by volatility modulated {L}\'evy-driven
  {V}olterra processes.
\newblock {\em Bernoulli\/}~{\em 19\/}(3), 803--845.

\bibitem[\protect\citeauthoryear{Barndorff-Nielsen, Maejima, and
  Sato}{Barndorff-Nielsen et~al.}{2006}]{idprocessesmaejimasatoole}
Barndorff-Nielsen, O.~E., M.~Maejima, and K.~Sato (2006).
\newblock Infinite divisibility for stochastic processes and time change.
\newblock {\em Journal of Theoretical Probability\/}~{\em 19\/}(2), 411--446.

\bibitem[\protect\citeauthoryear{Basse-O'Connor}{Basse-O'Connor}{2013}]{gammaandreas}
Basse-O'Connor, A. (2013).
\newblock Some properties of a class of continuous time moving average
  processes.
\newblock {\em Proceedings of the 18th EYSM\/}, 59--64.

\bibitem[\protect\citeauthoryear{Basse-O'Connor, Nielsen, Pedersen, and
  Rohde}{Basse-O'Connor et~al.}{2017}]{BasseNielPedRoh17}
Basse-O'Connor, A., M.~Nielsen, J.~Pedersen, and V.~Rohde (2017).
\newblock A continuous-time framework for arma processes.
\newblock {\em ArXiv e-prints\/}.

\bibitem[\protect\citeauthoryear{Basse-O'Connor and
  Rosi{\'n}ski}{Basse-O'Connor and Rosi{\'n}ski}{2013}]{BasseRosinski13}
Basse-O'Connor, A. and J.~Rosi{\'n}ski (2013).
\newblock Characterization of the finite variation property for a class of
  stationary increment infinitely divisible processes.
\newblock {\em Stochastic Processes and their Applications\/}~{\em 123\/}(6),
  1871--1890.

\bibitem[\protect\citeauthoryear{Brockwell and Davis}{Brockwell and
  Davis}{1986}]{Brockwell:1986:TST:17326}
Brockwell, P.~J. and R.~A. Davis (1986).
\newblock {\em Time Series: {T}heory and Methods}.
\newblock Springer-Verlag New York, Inc.

\bibitem[\protect\citeauthoryear{Brockwell and Lindner}{Brockwell and
  Lindner}{2009}]{Brockwell20092660}
Brockwell, P.~J. and A.~Lindner (2009).
\newblock Existence and uniqueness of stationary {L}\'evy-driven {CARMA}
  processes.
\newblock {\em Stochastic Processes and their Applications\/}~{\em 119\/}(8),
  2660 -- 2681.

\bibitem[\protect\citeauthoryear{Cohen and Maejima}{Cohen and
  Maejima}{2011}]{CohMaej2011}
Cohen, S. and M.~Maejima (2011).
\newblock Selfdecomposability of moving average fractional {L}\'evy processes.
\newblock {\em Statistics and Probability Letters\/}~{\em 81\/}(11),
  1664--1669.

\bibitem[\protect\citeauthoryear{Comte and Renault}{Comte and
  Renault}{1996}]{invnoncausacomte}
Comte, F. and E.~Renault (1996).
\newblock Noncausality in continuous time models.
\newblock {\em Econometric Theory\/}~{\em 12}, 215--256.

\bibitem[\protect\citeauthoryear{Duistermaat and Kolk}{Duistermaat and
  Kolk}{2010}]{DuisKolk10}
Duistermaat, J. and J.~Kolk (2010).
\newblock {\em Distributions: Theory and Applications}.
\newblock Cornerstones. Birkh{\"a}user Boston.

\bibitem[\protect\citeauthoryear{Ferrazzano and Fuchs}{Ferrazzano and
  Fuchs}{2013}]{ferrazzano2013}
Ferrazzano, V. and F.~Fuchs (2013).
\newblock Noise recovery for {L}\'evy-driven {CARMA} processes and
  high-frequency behaviour of approximating {R}iemann sums.
\newblock {\em Electronic Journal of Statistics\/}~{\em 7}, 533--561.

\bibitem[\protect\citeauthoryear{Kaminska}{Kaminska}{1997}]{Kaminska1997}
Kaminska, A. (1997).
\newblock On {M}usielak-{O}rlicz spaces isometric to {$L_2$} or {$L_{\infty}$}.
\newblock {\em Collectanea Mathematica\/}~{\em 48\/}(4-5-6), 563--569.

\bibitem[\protect\citeauthoryear{Pedersen and Sauri}{Pedersen and
  Sauri}{2015}]{PedSau15}
Pedersen, J. and O.~Sauri (2015).
\newblock On {L}\'evy semistationary processes with a gamma kernel.
\newblock In {\em XI Symposium on Probability and Stochastic Processes},
  Volume~69 of {\em Progress in Probability}, pp.\  217--239. Springer
  International Publishing.

\bibitem[\protect\citeauthoryear{Rajput and Rosi{\'n}ski}{Rajput and
  Rosi{\'n}ski}{1989}]{rajputrosinski}
Rajput, B.~S. and J.~Rosi{\'n}ski (1989).
\newblock Spectral representations of infinitely divisible processes.
\newblock {\em Probability Theory and Related Fields\/}~{\em 82\/}(3),
  451--487.

\bibitem[\protect\citeauthoryear{Rao and Ren}{Rao and Ren}{1994}]{Orliczrao}
Rao, M.~M. and Z.~D. Ren (1994).
\newblock {\em Theory of {O}rlicz spaces}.
\newblock New York: M. Dekker.

\bibitem[\protect\citeauthoryear{Sato}{Sato}{2006}]{satostochinadditive}
Sato, K. (2006).
\newblock Additive processes and stochastic integrals.
\newblock {\em Illinois J. Math\/}~{\em 50\/}(1 - 4), 825 -- 851.

\bibitem[\protect\citeauthoryear{Thuong}{Thuong}{2000}]{troungvan}
Thuong, T.~V. (2000).
\newblock Some colletions of functions dense in an {O}rlicz space.
\newblock {\em Acta Mathematica Vietnamica\/}~{\em 25\/}(2), 195--208.

\end{thebibliography}

\end{document}